\documentclass[11pt,amsfonts, epsfig]{amsart}

\usepackage{amsmath, amscd, amssymb}
\usepackage[frame,cmtip,arrow,matrix,line,graph,curve]{xy}
\usepackage{epsfig}
\usepackage{graphpap, color}
\usepackage[mathscr]{eucal}
\usepackage{mathrsfs}

\usepackage{pstricks}
\usepackage{color}
\usepackage{cancel}

\newcommand{\ds}{\displaystyle}

\numberwithin{equation}{section}

\newcommand{\codim}{\operatorname{codim}}

\def\sK{{\mathscr K}}
\def\sT{{\mathscr T}}
\def\sO{{\mathscr O}}

\def\sQ{{\mathscr Q}}

\def\sO{\mathscr{O}}

\def\sE{\mathscr{E}}

\def\sG{\mathscr{G}}

\newcommand{\PP}{\mathbb{P}}

\newcommand{\ZZ}{\mathbb{Z}}
\newcommand{\proj}{\mathbb P}

\newcommand{\bk}{\mathbf{k}}

\newcommand{\kk}{\bk}

\newcommand{\cal}{\mathcal}

\def\cF{{\cal F}}

\def\cQ{{\cal Q}}

\def\cT{{\cal T}}

\newcommand{\GL}{{\rm GL}}

\def\lra{\longrightarrow}

\newcommand{\ee}{{\bf e}}

\def\begeq{\begin{equation}}
\def\endeq{\end{equation}}
\def\and{\quad{\rm and}\quad}

\def\sub{\subset}

\def\and{\quad\text{and}\quad}

 \DeclareMathOperator{\Hom}{Hom}
\DeclareMathOperator{\HHom}{\mathcal{H}\it{om}}

 \DeclareMathOperator{\rank}{rank}

\let\lab=\label

\newtheorem{prop}{Proposition}[section]
\newtheorem{theo}[prop]{Theorem}
\newtheorem{lemm}[prop]{Lemma}
\newtheorem{coro}[prop]{Corollary}

\newtheorem{defi}[prop]{Definition}

\theoremstyle{definition}
\newtheorem{say}[prop]{}

\DeclareMathOperator{\EExt}{\mathcal{E}{\it xt}}

\def\Po{{\mathbb P^1}}

\def\PP{{\mathbb P}}

\let\lab=\label

\def\sO{{\mathscr O}}

\def\lab#1{\label{#1}[{#1}]\  }

\def\lab{\label} 

\def\beq{\begin{equation}}
\def\eeq{\end{equation}}

\title[A Compatification of the Space of Algebraic Maps]
{A Compactification of the Space of Algebraic Maps from
$\mathbb{P}^1$ to $\mathbb{P}^n$}

\author{Yi Hu}
\address{Department of Mathematics, University of Arizona, Tucson, AZ 85721 USA.}
\email{yhu@math.arizona.edu}
\author{Jiayuan Lin}
\address{Department of Mathematics, SUNY at Canton, NY 13617, USA}
\email{linj@canton.edu}
 \author{Yijun Shao}
\address{ Department of Mathematics, University of Arizona, Tucson,  AZ 85721 USA}
\email{yshao@math.arizona.edu}
\date{}

\begin{document}
\maketitle

\begin{abstract}
We provide a natural smooth projective compactification of the
space of algebraic maps from $\proj^1$ to $\proj^n$ by adding a
divisor with simple normal crossings.
\end{abstract}

\section{Introduction}

Fix a vector space $V$ of dimension $n+1$.  Let $N_d$ be the Quot
scheme parameterizing the exact sequences
$$0 \lra \sO_{\Po}(-d) \stackrel{f}{\lra} V\otimes \sO_{\Po} \lra
\cQ \lra 0,$$  where $\cQ$ is a coherent sheaf over $\Po$ of degree
$d$ and rank $n$. The locus of points of $N_d$ where $\cQ$ is
locally free can be identified with the space $\mathring{N}_d$ of
algebraic maps of degree $d$ from $\Po$ to $\proj (V)$. The
boundary $N_d \setminus \mathring{N}_d,$  consisting of points
parameterizing the exact sequences above where $\cQ$
is not locally free, has rather complicated singularities. One of
the main theme of the current paper is to resolve the
singularities of $N_d \setminus \mathring{N}_d$. For this, we find
that $N_d \setminus \mathring{N}_d$ comes equipped with a natural
filtration by subschemes $$  Z_{d,0}\subset Z_{d,1}\subset \cdots
\subset Z _{d,d-1}=N_d\setminus \mathring{N}_d $$  where $Z_{d,k}$
is supported on the subset $$ \{[f] \in N_d \;|\; \text{the
torsion part of }\cQ_{\rm}\text{ has degree}\geq d-k\} $$ for all
$0 \le k \le d-1$. As noted in pages 4738-39 of \cite{Hu2003}, it
is expected that one can successively blow up $N_d$ along the
subschemes $Z_{d,0}, Z_{d,1},  \cdots, Z _{d,d-1}$ such that
 the resulting final scheme $M_d$ is smooth and the boundary
 $M_d \setminus \mathring{N}_d $ is a divisor with normal crossings.
 The main purpose of the current paper is to execute this
conjectural construction in details.

We now briefly outline the main constructions. For each integer
$m\geq 0$, we let $$\rho_{f,m}: \Hom(V,\sO_{\Po}(m))\lra
\Hom(\sO_{\Po}(-d),\sO_{\Po}(m))$$ be the homomorphism obtained by
applying $\Hom(-,\sO_{\Po}(m))$ to $f$.
Fix any $m \ge d-1$. We set
$$W_m:=\Hom(V, \sO_{\Po}(m)) \and V_{d+m}:=H^0(\sO_{\Po}(d+m)).$$
Let $0 \le k \le d-1$. Then the exterior power
$\bigwedge^{k+2+m}\rho_{d,m}$ is a section of
\[
\Hom(\bigwedge^{k+2+m}W_{m},\bigwedge^{k+2+m}V_{d+m})\otimes
\sO_{N_d}(k+2+m).
\]
We proved that the scheme of zeros of
$\bigwedge^{k+2+m}\rho_{d,m}$ is independent of $m \ge d-1$ and it
is by definition the subscheme  $Z_{d,k}$.

For any two vector spaces $E$ and $F$, we let $S(E, F)$ denote
$\proj(\Hom(E,F))$. Then our main theorem reads

\begin{theo}\lab{main:grand} The variety
$M_d$ is a  compactification of $\mathring{N}_d$ such that the
following hold.
\begin{enumerate}
\item $M_d$ is isomorphic to the closure of the graph of the
rational map
\[
\begin{split}
N_d & \dashrightarrow \prod_{l=0}^{d-1} S(\bigwedge^{m+2+l}
W_{m},\bigwedge^{m+2+l} V_{d+m})\\ [f] &\mapsto
([\bigwedge^{m+2}\rho_{f,m}],[\bigwedge^{m+3}\rho_{f,m}],\dots,[\bigwedge^{m+d+1}\rho_{f,m}])
\end{split}
\]
for all $m \ge d-1$;  \item $M_d$ is a nonsingular projective
variety; \item The complement $M_d \setminus \mathring{N}_d $ is a
divisor with simple normal crossings.
\end{enumerate}
\end{theo}

In the course of proofs, we discover that our space $M_d$
possesses structures strikingly similar to the classic and modern
theories on complete collineations, complete correlations, and
complete quadrics. Indeed, our proofs rely on Vainsencher's
construction of the spaces of  the complete collineations
\cite{Vainsencher84, Vainsencher82}. The beautiful stories on
these complete objects went all the way back to the works of
Schubert in 19th century, to the works of\footnote{At the risk of
inadvertently omitting many authors who made important
contributions to these areas, we mention only a few.} Severi, Van
de Waerden, Semple, Tyrrell in the early and middle of the last
century, and to the modern treatments, refinements and advances of
Laskov \cite{Laksov85, Laksov88}, Vainsencher \cite{Vainsencher82,
Vainsencher84}, Thorup-Kleiman \cite{TK}, De Concini-Procesi
\cite{DP}, Demazure, and De Concini-Procesi-Goresky-MacPherson
\cite{DPGM} in the 1980's. Needless to say, this way of producing
good compactifications is nowadays very standard with the
Fulton-MacPherson compactification (\cite{FM}) and
Procesi-MacPherson compactification (\cite{MP}) being  the prime
examples. Related works in 1990's include the works of De
Concini-Procesi schools on hyperplane arrangements and the works
of Bifet-De Concini-Procesi on regular embeddings. Lately there
are further works in this direction by Wenchuan Hu and Li Li.

We believe that there are some lurking geometric objects,
analogous to the above classic complete objects, for which our
space $M_d$ is a parameter space. This is being pursued in a
forthcoming publication \cite{HS}.

Using Quot schemes of coherent sheaves of higher co-ranks, similar
spaces can be constructed to provide good compactifications of the
spaces of maps from the smooth rational curve to Grassmannians.
This has been  carried out in the third author's Ph.D dissertation
\cite{Shao}.

Two further problems to consider are to generalize to higher genus
curves and to compare our compactification with the Kontsevich
moduli space of stable maps (see, for example, \cite{Alexeev94,
BM, FP, KM}).

The current version of this work is resulted from a substantial
revision of the previous version. We would like to thank Gerd
Faltings for pointing out an error and the anonymous referee for pointing
out a gap in an early version and for his many helpful comments.
Many useful comments from
colleagues are also gratefully acknowledged. During the
preparation and the revision of this paper, the first author was
partially supported by NSA and NSF DMS 0901136.

 Throughout the
paper, we will work with a fixed algebraically closed base field
of characteristic zero, unless otherwise stated.

\section{Conventions and Terminology}

\begin{say}
From the introduction, $N_d$ is the Quote scheme parameterizing the exact sequences
\beq\lab{pointOfQuot} 0 \lra \sO_{\Po}(-d) \stackrel{f}{\lra}
V\otimes \sO_{\Po} \lra \cQ \lra 0\eeq where $\cQ$ is a coherent
sheaf over $\Po$ of degree $d$ and rank $n$.
We will
denote the corresponding point of \eqref{pointOfQuot} by
$$[f] \in N_d.$$
Note that we have the
 following identification
 $$N_d = \proj
(\Hom(\sO_{\Po}(-d),V\otimes \sO_{\Po})).$$ 
The Quote scheme $N_d$  comes equipped with
the universal family which is an exact sequence of coherent sheaves on $\mathbb{P}^1\times N_d$:
\begin{equation}\label{universalFamily}
0\lra \sO_{\Po}(-d)\otimes\sO_{N_d}(-1)\lra V\otimes \sO_{\Po\times N_d}\lra \sQ \lra 0.
\end{equation}
where $\sQ$ is a coherent sheaf of rank $n$, of relative degree $d$ and flat over $N_d$.
The restriction of \eqref{universalFamily} to the fiber of the projection $\Po\times
N_d\to N_d$ at the point $[f]$ is exactly the exact sequence (\ref{pointOfQuot}).
\end{say}

\begin{say}
Alternatively,  we may realize $N_d$ as
$$N_d=\{[f_0,\cdots ,f_n] \;|\; f_0,\cdots ,f_n  \in H^0(\sO_{\Po}(d)), \;\hbox{not all are zero\;}\}.$$
From this perspective, the boundary strata of $N_d \setminus \mathring{N}_d$ are
 $$C_{d,k}=\{[f_0,\cdots ,f_n]  \;|\;
f_0,\cdots ,f_n \;\text{have a common factor of degree}\geq d-k\}
$$ for all $0 \le k \le d-1$. More details along this line is to
be given in \S \ref{section:commonfactor}.
\end{say}

\begin{say}  For any non-negative integer $m$, we set
$$V_m:=H^0(\sO_{\Po}(m)),$$
$$W_m:=\Hom(V, \sO_{\Po}(m)).$$
In addition, we have the identifications
$$V_{d+k}=\Hom(\sO_{\Po}(-d),\sO_{\Po}(k)),$$
$$N_d = \mathbb{P}(V_d\otimes V).$$
\end{say}

\section{Resultant Homomorphisms}\lab{resultant}

\subsection{Resultant homomorphism and degree of
torsion}\lab{resul-torsion}

\begin{say} Consider the exact sequence \eqref{pointOfQuot}
$$0 \lra \sO_{\Po}(-d) \stackrel{f}{\lra} V\otimes \sO_{\Po} \lra
\cQ \lra 0.$$ For each integer $k\geq 0$,  recall from the introduction that
\beq\lab{resultHom} \rho_{f,k}: \Hom(V,\sO_{\Po}(k))\lra
\Hom(\sO_{\Po}(-d),\sO_{\Po}(k))\eeq is the map obtained by
applying $\Hom(-,\sO_{\Po}(k))$ to $f$. We call $\rho_{f,k}$ the
{\it $k$th resultant homomorphism} of $f$.
\end{say}

\begin{prop}\lab{rank-cT} Let  $\cT$ be the
torsion submodule of $\cQ$. Then we have
\[ \rank \rho_{f,k}\left\{\begin{array}{ll}
=k+1+d-\deg \cT, &\text{if }\ k\geq d-\deg \cT-1,\\
\geq 2(k+1),&\text{if }\ 0\leq k\leq d-\deg \cT-1.
\end{array}\right.\]
In particular, $\rank \rho_{f,k} = 2(k+1)$ when $k=  d-\deg
\cT-1$.
\end{prop}
\begin{proof}
First, we write $\cQ=\cF\oplus \cT$ such that $\cF$ is a locally
free sheaf of rank $n$ and degree $(d-\deg \cT)$. By applying the
functor $\Hom(-,\sO_{\Po}(k))$ to the exact sequence
\[
0\lra \sO_{\Po}(-d)\lra V\otimes \sO_{\Po}\lra \cQ\lra 0,
\]
we get an exact sequence
\[
0\lra \Hom(\cQ,\sO_{\Po}(k))\lra
\Hom(V,\sO_{\Po}(k))\stackrel{\rho_{f,k}}{\lra}
\Hom(\sO_{\Po}(-d),\sO_{\Po}(k)).
\]
Thus we have
\[
\begin{split}
\rank \rho_{f,k}&=\dim \Hom(V,\sO_{\Po}(k))-\dim
\Hom(\cQ,\sO_{\Po}(k))\\
&=\dim H^0(V^\vee(k))-\dim \Hom(\cF,\sO_{\Po}(k))\\
&=(n+1)(k+1)-\dim H^0(\cF^\vee(k)).
\end{split}
\]
We write $\cF$ as $\bigoplus_{i=1}^n\sO_{\Po}(d_i)$ with $d_i\geq
0$ and $\sum_{i=1}^n d_i=d-\deg \cT$. Then
\[
H^0(\cF^\vee(k))=H^0(\bigoplus_{i=1}^n\sO_{\Po}(k-d_i))=\bigoplus_{i=1}^n
H^0(\sO_{\Po}(k-d_i)).
\]
Consequently,
\[
h^0(\cF^\vee(k))=\sum_{i=1}^n\max\{k-d_i+1,0\}.
\]
If $k\geq d-\deg \cT-1$, then $k-d_i+1\geq 0$ for all $i$. In this
case
\[
h^0(\cF^\vee(k))=\sum_{i=1}^n(k-d_i+1)=n(k+1)-(d-\deg \cT)
\]
and it implies
\[
\rank \rho_{f,k}=k+1+d-\deg \cT.
\]
This proves the first case of the proposition.

On the other hand, if $0\leq k\leq d-\deg \cT-1$, then we claim
that \[ h^0(\cF^\vee(k))\le (n-1)(k+1).
\]
There are two cases to consider. First, $k-d_i+1 \ge 0$ for all
$i$. In this case,
\[
h^0(\cF^\vee(k))= \sum_{i=1}^n(k-d_i+1)=n(k+1)-(d-\deg \cT)\leq
(n-1)(k+1).
\]
Secondly, $k-d_i+1 < 0$ for some $i$, say for $i=1$. In this case,
\[
h^0(\cF^\vee(k))= \sum_{i=2}^n \max\{k-d_i+1,0\} \le \sum_{i=2}^n
(k+1) =(n-1)(k+1).
\]
Thus, in either case,
 $\rank
\rho_{f,k}\geq(n+1)(k+1)-(n-1)(k+1)=2(k+1)$.

This completes the proof.
\end{proof}


\begin{coro}\lab{le2k+1}
Let the notations be as in above.
\begin{enumerate}
\item  Assume that $\rank \rho_{f,k} \le 2k+1$.  Then we have
$$\deg \cT \ge d-k \and \rank \rho_{f,k}=k+1+d-\deg \cT.$$
Further,  for all $m \ge k-1$, $$\rank \rho_{f,m}-\rank
\rho_{f,k}=m-k;$$ \item For  any fixed $l \ge 0$, $$\deg \cT=d-k
\;\hbox{if and only if}\; \rank\rho_{f,k+l}=2k+1+l;$$ \item For
any fixed $l \ge 0$, $$\deg \cT \geq d-k \;\hbox{if and only if}\;
\rank\rho_{f,k+l} \leq 2k+1+l.$$
\end{enumerate}
\end{coro}
\begin{proof} (1). By Proposition \ref{rank-cT}, we must have $k > d-\deg
\cT-1$, that is $\deg \cT \ge d-k$. In this case, $\rank
\rho_{f,k}=k+1+d-\deg T$. Further, when $m \ge k-1$, because
$k+1+d-\deg \cT = \rank \rho_{f,k} \le 2k+1$, we must also have $m
\ge d -\deg \cT -1$. Hence by Proposition \ref{rank-cT} again,
$$\rho_{f,m}-\rho_{f,k}=(m+1+d-\deg \cT)-(k+1+d-\deg \cT)=m-k. $$

(2). $\deg \cT=d-k$ if and only if (by Proposition \ref{rank-cT})
$\rank \rho_{f,k}=2k+1$ if and only if (by (1)) $\rank
\rho_{f,k+l}=2k+1 +l$.

(3). $\deg \cT \ge d-k$ if and only if (by Proposition
\ref{rank-cT}) $\rank \rho_{f,k}\le 2k+1$ if and only if(by (1))
$\rank \rho_{f,k+l} \le 2k+1 +l$.
\end{proof}

\subsection{Resultant homomorphism and degree of common
factor}\lab{section:commonfactor}

\begin{say}
In this subsection, we reinterpret the results of \S
\ref{resul-torsion} in terms of $``$common factors of
polynomials$"$. For this, we use the following identification
$$\Hom(\sO_{\Po}(-d),V)=V\otimes\Hom(\sO_{\Po}(-d),\sO_{\Po})=V\otimes
H^0(\sO_{\Po}(d)).$$ Then we let $\{x,y\}$ be a basis for
$H^0(\sO_{\Po}(1))$ and let $\{e_0,\dots,e_n\}$ be a basis for
$V$. This way,  any $f\in \Hom(\sO_{\Po}(-d),V)$ can be written as
\[
f=e_0\otimes f_0+\dots+ e_n\otimes f_n
\]
where $f_i=\sum_{j=0}^d a_{ij}x^{d-j}y^j$ are homogeneous
polynomials in $x,y$ of degree $d$. Recall that under this view of
points of the space $N_d$, $\deg \cT$ is simply the degree of the
greatest common factors of $f_0, \ldots, f_n$.
\end{say}

\begin{say}\lab{bases} Consider the $k$th resultant homomorphism \eqref{resultHom}
 $$\rho_{f,k}: \Hom(V,\sO_{\Po}(k))\lra
\Hom(\sO_{\Po}(-d),\sO_{\Po}(k)).$$ We abbreviate $\rho_{f,k}$ as
$\rho_{f,k}: W_k \lra V_{d+k}$ where
\[
\begin{split} & W_k:=\Hom(V,\sO_{\Po}(k))= V^\vee \otimes H^0(\sO_{\Po}(k))
\\
&
V_{d+k}:=\Hom(\sO_{\Po}(-d),\sO_{\Po}(k))=H^0(\sO_{\Po}(d+k)).\end{split}\]
 Using the basis $\{e_i^\vee\otimes x^ky^{k-j}
:i=0,\dots,n;\ j=0,\dots,k\}$ for $W_k$ and the basis
$\{x^{d+k-j}y^j: j=0,\dots,d+k\}$ for $V_{d+k}$, the linear map
$\rho_{f,k}$ is represented by the following
$(k+1)(n+1)\times(d+k+1)$ matrix \beq\lab{matrixA}
A_{f, k}=\begin{pmatrix} a_{00} & a_{01} & \cdots & a_{0d} \\
\vdots & \vdots & \cdots & \vdots \\
a_{n0} & a_{n1} & \cdots & a_{nd} \\
& a_{00} & a_{01} & \cdots & a_{0d} \\
& \vdots & \vdots & \cdots & \vdots \\
& a_{n0} & a_{n1} & \cdots & a_{nd} \\
& & \cdots & \cdots & \cdots & \cdots \\
& & & \cdots & \cdots & \cdots & \cdots \\
& & & & a_{00} & a_{01} & \cdots & a_{0d}  \\
& & & & \vdots & \vdots & \cdots & \vdots \\
& & & & a_{n0} & a_{n1} & \cdots & a_{nd}
\end{pmatrix}
\eeq which acts on elements of  $W_k$ by multiplication from the
right.
\end{say}

Corollary \ref{le2k+1} takes the following form in this setting.

\begin{coro}\lab{le2k+1:commonFactor} Let $g$ be a greatest common factor of
$f_0, \cdots, f_n$.
\begin{enumerate}
\item Assume that $\rank A_{f,k} \le 2k+1$. Then we have
$$\deg g \ge d-k \and \rank A_{f,k}=k+1+d-\deg g.$$  Further,  for all $m
\ge k-1$,
$$\rank A_{f,m}-\rank A_{f,k}=m-k;$$ \item For  any fixed $l \ge
0$, $$\deg g=d-k \;\hbox{ if and only if } \;\rank
A_{f,k+l}=2k+1+l;$$  \item For  any fixed $l \ge 0$,  $$\deg g
\geq d-k \;\hbox{ if and only if } \;\rank A_{f,k+l} \leq
2k+1+l.$$
\end{enumerate}
\end{coro}

We remark here that Proposition \ref{rank-cT} implies Kakie's
Proposition 3, \cite{Kakie}.

\section{Determinantal Subschemes}\label{zeroLoci}

\subsection{The universal resultant homomorphisms}

\begin{say} Let $\pi: \Po\times N_d\to N_d$ denote the second projection.
For each integer $m \ge 0$,
by applying the functor $\pi_*\HHom(-,\sO_{\Po}(m))$ to the homomorphism
$\sO_{\Po}(-d)\otimes\sO_{N_d}(-1)\to V\otimes \sO_{\Po\times N_d}$ (which comes from the universal
family (\ref{universalFamily})),  we obtain { a nowhere zero $\sO_{N_d}$-homomorphism
\[
\rho_{d,m}: W_m\to V_{d+m}\otimes \sO_{N_d}(1).
\]
It is a nowhere zero section of $\Hom(W_m,V_{d+m})\otimes \sO_{N_d}(1)$ whose restriction
to every point $[f]\in N_d$ is $\rho_{f,m}$.} We call $\rho_{d,m}$ the $m$-th {\it universal
resultant homomorphism}.
\end{say}

\begin{say}\lab{Zdkm}
 For $m, k \geq 0$,  the exterior power
$\bigwedge^{k+2+m}\rho_{d,m}$ is a section of
\[
\Hom(\bigwedge^{k+2+m}W_{m},\bigwedge^{k+2+m}V_{d+m})\otimes
\sO_{N_d}(k+2+m).
\]
Fix any $1 \le k \le d-1$. Then using Corollary \ref{le2k+1}, one
checks that the scheme $Z_{d,k;m}$ of zeros of
$\bigwedge^{k+2+m}\rho_{d,m}$ is supported on $C_{d,k}$ whenever
$m \ge k$. The ideal sheaf  $I_{d,k;m}$ of $Z_{d,k;m}$ is the
image of the induced homomorphism
\[
\Hom(\bigwedge^{k+2+m}W_{m},\bigwedge^{k+2+m}V_{d+m})^\vee\otimes
\sO_{N_d}(-k-2-m)\twoheadrightarrow I_{d,k;m}\subset \sO_{N_d}.
\]
\end{say}

\begin{say} We suspect that for any fixed $1 \le k \le d-1$,
\beq I_{d,k;m}= I_{d,k;k}\eeq holds whenever $m \ge k$. This would
imply that  $I_{d,k;m}$ with $m \ge k$ all endow the same scheme
structure on $C_{d,k}$.  Rather than proving this, we will show
the weaker Proposition \ref{weaker-sameIdeals} below, which
already suffice for our purpose. To pave the way for its proof, we
need some preparation.
\end{say}

\begin{say}\lab{minorIdeals}
Let $R$ be a ring and $A$ a $p \times q$ matrix over $R$. We let
$I_l(A)$ be the ideal generated by all $l \times l$ minors of $A$
with $1 \le l \le p, q$.  Suppose $B$ is an invertible $p \times
p$ matrix and $C$ is an invertible $q \times q$ matrix. Then one
checks directly that \beq\lab{sameMinorIdeals} I_l(A) = I_l(BA) =
I_l(AC)\eeq for all $1 \le l \le p, q$. In more concrete terms,
\eqref{sameMinorIdeals} means that the following three operations
on the matrix $A$ preserve the ideal $I_l(A)$:
\begin{enumerate}
\item multiply a row or a column by units; \item interchanging two
rows or two columns; \item multiply one row (column) by an element
of $R$ and add the result to another row (column).
\end{enumerate}
\end{say}

\begin{prop}\lab{weaker-sameIdeals} Fix $1\leq k\leq d-1$. Then
\begin{enumerate}
\item  $I_{d,0;m}=I_{d,0;0}$ for all $m\geq 0$; \item
$I_{d,k;m}=I_{d,k;d-1}$ for all $m \geq d-1$.
\end{enumerate}
\end{prop}
\begin{proof}
 Consider any nonzero $f\in \Hom(\sO_{\Po}(-d),V)$. Recall that using the bases as
chosen in \ref{bases}, we express it as $f=e_0\otimes f_0+\dots+
e_n\otimes f_n $ where $f_i=\sum_{j=0}^d a_{ij}x^{d-j}y^j$. This
way, the coefficients $(a_{ij})$  become the homogeneous
coordinates of $N_d$. Observe in addition  that for any fixed $m
\ge 0$, $(a_{ij})$ are also the entries in the first block of the
matrix $A_{f,m}$ of the $m$-version of \eqref{matrixA}. We regard
$A_{f,m}$ as a matrix over the polynomial ring $\kk[a_{ij}]$.
Observe that the ideal sheaf $I_{d,k;m}$ coincides with the sheaf
 $(I_{k+2+m}(A_{f,m}))^{\sim}$ associated to the module $I_{k+2+m}(A_{f,m})$. Our strategy of the
proof is to cover $N_d$ by the standard affine open subsets and
prove the statements over the open subsets.

We first localize to the affine open set $U_0=(a_{00}\neq 0)$. By
using the affine coordinates $b_{ij}=a_{ij}/a_{00}$, the matrix
$A_{f,m}$ is reduced to $B_{f,m}$ such that all the entries
$a_{ij}$ are replaced by $b_{ij}$ except that $a_{00}$ is replaced
by 1. Then the localization of each ideal $I_{k+2+m}(A_{f,m})$ to
$U_0$ is $I_{k+2+m}(B_{f, m})$. Now using $b_{00}=1$, we can
eliminate the entries $b_{i0}$ by some appropriate row operations
and reduce the matrix $B_{f, m}$ to
\[
C_{f, m}=\begin{pmatrix} 1 & c_{01} & \cdots & c_{0d} \\
0 & c_{11} & \cdots & c_{1d} \\
\vdots & \vdots & \cdots & \vdots \\
0 & c_{n1} & \cdots & c_{nd} \\
& 1 & c_{01} & \cdots & c_{0d} \\
& 0 & c_{11} & \cdots & c_{1d} \\
& \vdots & \vdots & \cdots & \vdots \\
& 0 & c_{n1} & \cdots & c_{nd} \\
& & \cdots & \cdots & \cdots & \cdots \\
& & & \cdots & \cdots & \cdots & \cdots \\
& & & & 1 & c_{01} & \cdots & c_{0d}  \\
& & & & 0 & c_{11} & \cdots & c_{1d} \\
& & & & \vdots & \vdots & \cdots & \vdots \\
& & & & 0 & c_{n1} & \cdots & c_{nd} \\
\end{pmatrix}
\]
where $c_{0j}=b_{0j}, 1 \le i \le d$ and
$c_{ij}=b_{ij}-b_{i0}b_{0j}, 1 \le i \le n, 1 \le j \le d$.

We are now ready to prove (1) over the open subset $U_0$. One sees
by direct calculations that
 $$I_2(C_{f,0})=\langle c_{ij} \;|\; 1 \le i \le n, 1 \le
j \le d  \rangle.$$ For the ideal $I_{2+m}(C_{f,m})$ with
$m\geq1$,  it is trivial that $$I_{2+m}(C_{f,m})\subset
I_2(C_{f,0}).$$ On the other hand, observe that $C_{f,m}$ has a
$(n+1+m)\times(d+m+1)$ submatrix of the form
\[
\begin{pmatrix} 1 & c_{01} & \cdots & c_{0d} \\
& 1 & c_{01} & \cdots & c_{0d} \\
&   & \ddots & \ddots & \ddots & \ddots \\
&   &  & 1 & c_{01} & \cdots & c_{0d} \\
&   &  & 0 & c_{11} & \cdots & c_{1d} \\
&   &  &\vdots & \vdots & \cdots & \vdots \\
&   &  & 0 & c_{n1} & \cdots & c_{nd} \\
\end{pmatrix}.
\]
Using the $(m+1)$ 1's on the diagonal, one easily finds
$(m+2)\times (m+2)$ minors such that their determinants are
$c_{ij}$ for all $1 \le i \le n$ and $1 \le j \le d$. This implies
that $I_2(C_{f,0})\subset I_{2+m}(C_{f,m})$. Thus,
$$I_{2+m}(C_m)=I_2(C_0)$$ for all $m\geq 1$. This proves (1) over
the open subset $U_0$.

We now turn to the statement (2) over $U_0$. Consider the matrix
$C_{f,m}$ with $m \ge d$. 
It is routine to check that the following holds: \beq \hbox{row}_i
+ \sum_{j=1}^d (c_{ij} \hbox{row}_{j(n+1)+i-1} + c_{0j}
\hbox{row}_{j(n+1)+i}) =0\eeq for all $2 \le i \le n+1$. This
means that we can eliminate $\hbox{row}_i$ for all $2 \le i \le
n+1$. We can also easily eliminate the entries $c_{0i}$ of the
first row by using the first column. This implies that
$I_{k+2+m}(C_{f,m})=I_{k+1+m}(C_{f,m-1})$. This process can be
repeated until we reach $I_{k+1+d}(C_{f, d-1})$. Thus we obtain
$$I_{k+2+m}(C_{f,m})= I_{k+1+d}(C_{f, d-1})$$
for all $m \ge d-1$. This completes the proof of (1) and (2) over
the open subset $U_0$.

To investigate (1) and (2) over the rest of open charts of $N_d$,
we use the symmetry of $N_d$. The group $$\GL(H^0(\sO_{\Po}(1)))
\times \GL(V)\cong \GL_2 \times \GL_{n+1}$$ acts on $N_d$. If
$g\in \GL(H^0(\sO_{\Po}(1)))$, it corresponds to  a change of
basis of $H^0(\sO_{\Po}(1))$; since it induces bases changes in
both $W_m$ and $V_{d+m}$, we see that $g$ acts on the matrix
$A_{f,m}$ by multiplying invertible matrices from both the left
and the right. Likewise, an element $g \in \GL(V)$ acts on
$A_{f,m}$ by multiplying an invertible matrix from the left. By
\ref{minorIdeals}, $g^*(I_l(A_{f,m}))=I_l(A_{f,m})$ for any $g\in
\GL(H^0(\sO_{\Po}(1))) \times \GL(V)$.

Now, for $[f] \in N_d \setminus U_0$, we have $a_{00}=0$. If one
of $a_{0i}$ is not zero, say $a_{0j} \ne 0$. It is routine to find
a basis change of $H^0(\sO_{\Po}(1))$ such that under the new
basis $a'_{00}\ne 0$. This means that there is $g \in
\GL(H^0(\sO_{\Po}(1)))$ such that $g\cdot [f] \in U_0$. If all of
$a_{0i}$ are zero, then there is $i \ge 1$ and $j$ such that
$a_{ij} \ne 0$. Then let $g \in \GL(V)$ correspond to
interchanging $e_0$ and $e_j$, we see that $g\cdot [f]$ places us
in the previous situation. In either case, by the invariance of
the ideals $g^*(I_l(A_{f,m}))=I_l(A_{f,m})$, we conclude that the
statements (1) and (2) hold everywhere in $N_d$.

This completes the proof.
\end{proof}

\subsection{Determinantal subschemes and their basic properties}

\begin{say} Let $1 \le k \le d-1$.
 We set
 $$I_{d,0}= I_{d,0;m}, \; m \ge 0 \and I_{d,k}=I_{d,k;m}, \; m \ge
 d-1.$$
 By Proposition \ref{weaker-sameIdeals}, these are well-defined.
\end{say}

\begin{defi}\lab{defn:zdk} For any $0 \le k \le d-1$,
we let $Z_{d,k}$ be the subscheme of $N_d$ defined by the ideal
$I_{d,k}$.
\end{defi}

By \ref{Zdkm}, $Z_{d,k}$ is supported on $C_{d,k}$.

\begin{say} Since $$N_d = \proj
(\Hom(\sO_{\Po}(-d),V\otimes \sO_{\Po})) \and N_k = \proj
(\Hom(\sO_{\Po}(-k),V\otimes \sO_{\Po})),$$ using the
identification $V_{d-k}=\Hom(\sO_{\Po}(-d),\sO_{\Po}(-k))$, we
obtain a natural morphism \beq\lab{phidk}
\begin{split}
\varphi_{d,k}\;: \;\;& \mathbb{P}(V_{d-k}) \times N_k \lra  N_d \\
 & ([h],[g]) \mapsto [g\circ h].
\end{split}
\eeq
 When $d=0$, $N_0=\proj(V)$. Also we have the identification
$N_d = \mathbb{P}(V_d\otimes V)$. A direct computation shows that
\end{say}

 \begin{prop} The morphism
$\varphi_{d,0}: \mathbb{P}(V_{d})\times \mathbb{P}(V) \lra
N_d=\mathbb{P}(V_d\otimes V)$ is the Segre embedding and its image
scheme is exactly $Z_{d,0}$.
 \end{prop}

In particular, this implies that $Z_{d,0}$ is smooth. For general
$\varphi_{d,k}$, we have

\begin{prop}\lab{prop:inverse}
The restriction of $\varphi_{d,k}$ to $\mathbb{P}(V_{d-k})\times
(N_k\setminus Z_{k,k-1})$ gives rise to an isomorphism
\beq\lab{inverse} \varphi_{d,k}': \mathbb{P}(V_{d-k})\times
(N_k\setminus Z_{k,k-1})\stackrel{\cong}{\lra} Z_{d,k}\setminus
Z_{d,k-1}.\eeq
\end{prop}
\begin{proof}  The idea of the proof is taken from the third author's thesis
\cite{Shao}. We give sufficient sketch here.

We just need  to produce the inverse to $\varphi_{d,k}'$.

First, recall that the Quot scheme $N_d$ comes equipped with a universal exact sequence of sheaves
over $\mathbb{P}^1\times N_d$ \beq\lab{universalNd} 0\to { \sE} \to V\otimes
\sO_{\mathbb{P}^1\times N_d}\to \sQ \to 0 \eeq where  { $\sE = \sO_{\mathbb{P}^1}(-d)\otimes
\sO_{N_d}(-1)$ and} $\sQ$ is a coherent sheaf of rank $n$, relative degree $d$ and flat over $N_d$.
Similarly, $\mathbb{P}(V_{d-k})$ comes equipped with a universal
exact sequence of sheaves over $\mathbb{P}^1\times
\mathbb{P}(V_{d-k})$ \beq\lab{universalPVd} 0\to
\sO_{\Po}(-d)\otimes
\sO_{\mathbb{P}(V_{d-k})}(-1)\to\sO_{\Po}(-k)\otimes\sO_{\mathbb{P}(V_{d-k})}
\to \sT \to 0 \eeq where $\sT$ is a torsion sheaf of relative
degree $d-k$ and is flat over $\mathbb{P}(V_{d-k})$.

Next, taking dual of the exact sequence \eqref{universalNd}, we
obtain \beq\lab{dual} V^\vee_{\mathbb{P}^1\times
N_d}\to\sO_{\Po}(d)\otimes \sO_{N_d}(1)\to \sG\to 0 \eeq where
$\sG=\EExt^1(\sQ,\sO_{\Po \times N_d})$.  Here and below, we use
$\EExt^1$ for $\EExt^1_{\sO_{\Po \times N_d}}$. Tensoring \eqref{dual}
by $\sO_{\Po}(m)$ for $m\gg 0$ and applying $\pi_*$ where $\pi:
\mathbb{P}^1\times N_d\to N_d$ is the projection map, we then
obtain \beq \pi_*(V^\vee_{\mathbb{P}^1\times N_d}(m)) \to
\pi_*(\sO_{\Po}(d+m)\otimes \sO_{N_d}(1)) \to \pi_*(\sG(m))\to 0
\eeq where the first map is simply
$$\rho_{d,m}: W_m = \pi_*(V^\vee_{\mathbb{P}^1\times N_d}(m)) \to
\pi_*(\sO_{\mathbb{P}^1}(d+m)\otimes \sO_{N_d}(1))=V_{d+m}\otimes
\sO_{N_d}(1).$$  Since $Z_{d,k}$ is the scheme of zeros of
$\bigwedge^{k+2+m}\rho_{d,m}$, we see that $\pi_*\sG(m)$ pulls
back to a locally free sheaf of rank $d-k$ over  $Z_{d,k}\setminus
Z_{d,k-1}$ for all $0 \le k \le d-1$. Set $Z_{d,-1}:=\emptyset$.
Then one checks directly that the disjoint union
$$ \bigsqcup_{k=0}^{d-1} (Z_{d,k}\setminus
Z_{d,k-1})\bigsqcup \;(N_d\setminus Z_{d,d-1})$$ is exactly the
flattening stratification of $\sG$ (cf. Lecture 8,
\cite{Mumford}).

Now, we let
$$\iota: \Po \times (Z_{d,k}\setminus Z_{d,k-1}) \lra \Po \times
N_d$$ be the inclusion. Then, the torsion sheaf $\iota^* \sG$ has
relative degree $d-k$ and is flat over $Z_{d,k}\setminus
Z_{d,k-1}$. Thus, by the universality of $\mathbb{P}(V_{d-k})$, we
obtain a morphism \beq \lab{mor1} Z_{d,k}\setminus Z_{d,k-1} \lra
\mathbb{P}(V_{d-k}).\eeq

What remains is to get a morphism from $Z_{d,k}\setminus Z_{d,k-1}$ to $N_k\setminus Z_{k,k-1}$.
For this, we pull back \eqref{universalNd} to $\Po \times (Z_{d,k}\setminus Z_{d,k-1})$. Since
$\sQ$ is flat over $N_d$, we get an exact sequence $$ 0\to \iota^*\sE \to V\otimes
\sO_{\mathbb{P}^1\times (Z_{d,k}\setminus Z_{d,k-1})}\to \iota^*\sQ \to 0. $$   Taking the dual to
the above sequence, we obtain
\begin{equation}\label{longSequence}
0\to (\iota^*\sQ)^\vee \to V^\vee\otimes \sO_{\mathbb{P}^1\times (Z_{d,k}\setminus Z_{d,k-1})} \to
(\iota^*\sE)^\vee \to { \EExt^1_Z(\iota^*\sQ,\sO)} \to 0,
\end{equation}
{  where $\EExt^1_Z:=\EExt^1_{\sO_{\mathbb{P}^1\times (Z_{d,k}\setminus Z_{d,k-1})}}$.
We have a canonical identification:
$$\EExt^1_Z(\iota^*\sQ,\sO)=\iota^*\EExt^1(\sQ,\sO)=\iota^*\sG.$$ To see this,
we  first dualize the universal exact sequence (\ref{universalNd}) to obtain
\[
0\to \sQ^\vee \to V^\vee\otimes \sO_{\mathbb{P}^1\times N_d} \to \sE^\vee \to \EExt^1(\sQ,\sO) \to
0.
\]
Then, we pull it back to $\mathbb{P}^1\times(Z_{d,k}\setminus Z_{d,k-1})$ to obtain
\[
V^\vee\otimes \sO_{\mathbb{P}^1\times (Z_{d,k}\setminus Z_{d,k-1})} \to \iota^*(\sE^\vee) \to
\iota^*\EExt^1(\sQ,\sO) \to 0.
\]
Since pulling-back and dualizing operations commute on locally free sheaves, we have a canonical
identification $(\iota^*\sE)^\vee=\iota^*(\sE^\vee)$ and hence a commutative diagram
\[
\xymatrix{ %
V^\vee\otimes\mathscr{O}_{\mathbb{P}^1\times (Z_{d,k}\setminus Z_{d,k-1})}\ar[r]\ar@{=}[d] &
\iota^*(\mathscr{E}^\vee)\ar[r]\ar@{=}[d] & \iota^*\mathscr{E}xt^1(\mathscr{Q},\mathscr{O})\ar[r]\ar@{:}[d] & 0 \\
V^\vee\otimes\mathscr{O}_{\mathbb{P}^1\times (Z_{d,k}\setminus Z_{d,k-1})}\ar[r] &
(\iota^*\mathscr{E})^\vee\ar[r] & \mathscr{E}xt^1_Z(\iota^*\mathscr{Q},\mathscr{O})\ar[r] & 0.
} %
\]
}
This gives the canonical identification in the third collumn, as desired.

 We now break up the sequence
(\ref{longSequence}) into two: \beq \lab{split1} 0\to (\iota^*\sQ)^\vee \to V^\vee\otimes
\sO_{\mathbb{P}^1\times (Z_{d,k}\setminus Z_{d,k-1})}\to \sK\to 0 \eeq and \beq\lab{split2} 0\to
\sK \to \iota^*\sE^\vee \to \iota^*\sG\to 0. \eeq Since $\iota^*\sG$ is flat over $Z_{d,k}\setminus
Z_{d,k-1}$, one checks that $\sK$ is locally free, hence  $(\iota^*\sQ)^\vee$ is locally free. Now
taking the dual of \eqref{split1}, we get an exact sequence of locally free sheaves
\[
0\to \sK^\vee\to V\otimes \mathcal{O}_{\mathbb{P}^1\times
(Z_{d,k}\setminus Z_{d,k-1})}\to (\iota^*\sQ)^{\vee\vee}\to 0.
\]
One easily calculates that 
$$\rank (\iota^*\sQ)^{\vee\vee}=n \and \deg
(\iota^*\sQ)^{\vee\vee}=k.$$
Thus, by the universality of the Quot scheme $N_k$, we obtain a
morphism
$$Z_{d,k}\setminus Z_{d,k-1}\to N_k\setminus Z_{k,k-1}.$$
Together with \eqref{mor1}, this gives rise to a morphism
$$\psi_{d,k}: Z_{d,k}\setminus
Z_{d,k-1}\to \mathbb{P}(V_{d-k})\times (N_k\setminus Z_{k,k-1})$$
which is  the inverse to $\varphi'_{d,k}$. This completes the
proof.
\end{proof}

As a consequence of this proposition, we see that
$Z_{d,k}\setminus Z_{d,k-1}$ is a locally closed smooth subvariety
of $N_d$ for all $k\geq 1$.

\section{Statements of Main Results and Their
Proofs}\lab{mainProofs}

\subsection{Describing the successive blowups}\lab{ourBlowups}

\begin{say} We use induction to
describe the iterated blowups of $N_d$  along the $Z_{d,0},
\ldots, Z_{d,d-1}$. We set $N_d^{-1}:=N_d$ and
$Z_{d,k}^{-1}:=Z_{d,k}$. For any $0\leq l\leq d-1$, we let $N_d^l$
be the blowup of $N_d^{l-1}$ along $Z_{d,l}^{l-1}$, $Z_{d,k}^l$
the proper transform of $Z_{d,k}^{l-1}$ when $k \ne l$, and
$Z_{d,l}^{l}$ the exceptional divisor of  $N_d^l \lra N_d^{l-1}$.
Observe that $Z_{d,k}^l$ is a divisor when $k \le l$ and
$Z_{d,k}^l$ is the blowup of $Z_{d,k}^{l-1}$ along $Z_{d,l}^{l-1}$
when $k >l$.
\end{say}

\begin{say} We denote the final blowup $N_d^{d-1}$ by $M_d$. We
aim to show that $M_d$ is smooth and the boundary $M_d \setminus
\mathring{N}_d$ is a divisor with normal crossings. The smoothness
of $M_d$ will follow by induction if the blowup center
$Z_{d,l+1}^{l}\sub N_d^l$ is smooth.  A technical key to prove
this is to show that the total transform in $N_d^l$ of
$Z_{d,l+1}^{l-1} \subset N_d^{l-1}$ is the scheme-theoretical
union of a Cartier divisor with the proper transform $Z_{d,l+1}^l
\subset N_d^l$, that is, locally we have
$$I_{Z_{d,l+1}} \cdot \sO_{N_d^l} = P \cdot I_{Z_{d,l+1}^l}$$
where $P$ is a principal ideal. Here by $I_Z$, we mean the ideal
sheaf of a subscheme $Z$. To this end, we will relate our blowups
to spaces of complete collineations and apply the related results
of Vainsencher \cite{Vainsencher84} and \cite{Vainsencher82}.
\end{say}

\subsection{Using the space of complete collineations}

\begin{say} Let $E$ and $F$ be vector spaces. Let $S(E,F)=\proj(\Hom(E,F))$
be the space of collineations from $E$ to $F$. It comes equipped
with a universal homomorphism
\[
u_{EF}: E\to F\otimes \sO_{S(E,F)}(1)
\]
For $k \ge 1$, set $D_k(E,F)$ to be the scheme of zeros of the
section
$$\bigwedge^{k+1}u_{EF}:  \Hom(\bigwedge^{k+1}E,\bigwedge^{k+1}F)\otimes
\sO_{S(E,F)}(k+1).$$
\end{say}

\begin{say}\lab{sEF} Let $r+1=\min\{\dim E, \dim
F\}$. Set $S^0:=S(E,F)$, $D_k^0:=D_k(E,F)$. We define the
following inductively.  For $1\leq l\leq r$, let $S^l$ be the blow
up of $S^{l-1}$ along $D_l^{l-1}$, $D_k^l$ the proper transform of
$D_{k}^l$ for $k\neq l$, and $D_l^l$ the exceptional divisor of
$S^l \lra S^{l-1}$.
\end{say}

\begin{say}
Although $D_k(E,F)$ is singular for $k \ge 2$, Vainsencher
\cite{Vainsencher84} shows that $D_l^{l-1}$ is smooth, thus $S^l$
is smooth. In particular, the final blowup space $S^r$ is smooth.
Further, he shows that $S^r$ parameterizes $``$complete
collineations$"$ from $E$ to $F$ (see \cite{Vainsencher84} for
more details). The property that we need from Vainsencher's
construction is the following
\end{say}

\begin{prop}\lab{vain} {\rm (Theorem 2.4 (8), \cite{Vainsencher84})}  Assume $1 \le l < k$.
\[
I_{D_k^{l-1}}\cdot\sO_{S^l}=I_{D_k^l}\cdot(I_{D_l^l})^{k-l+1}
\]
\end{prop}

\begin{say}\lab{relation} The relations between our blowups as described in \S \ref{ourBlowups}
and the spaces of $``$complete collineations$"$ are as follows.
Note that for any $m \ge 0$, the resultant homomorphism of
\eqref{resultHom}
$$
\rho_{f,m}: W_m \lra V_{d+m}$$ gives rise to an embedding
\[\begin{split}
N_d  & \lra   S(W_m,V_{d+m}) \\
 [f] & \mapsto [\rho_{f,m}].
 \end{split}
\]
Indeed, $N_d= S(W_0,V_d)$.  Further, $\rho_{d,m}$ is exactly the
pullback to $N_d$ of the universal map $u=u_{W_mV_{d+m}}$ on
$S(W_m,V_{d+m})$. Consequently, $\bigwedge^{l}u$ pulls back to
$\bigwedge^l\rho_{d,m}$ for all $l$.
Then it follows by definition that for all $0\leq k\leq d-1$ and
$m\geq d-1$
\[
N_d\cap D_{k+1+m}(W_{m},V_{d+m})=Z_{d,k}
\]
scheme-theoretically. Also it is easy to check that
$$\rank\rho_{f,m}\geq m+1$$
for all $m\geq 0$ and $[f]\in N_d$. Therefore, \beq\lab{empty}
N_d\cap D_1(W_m,V_{d+m})=\dots=N_d\cap D_m(W_m,V_{d+m})=\emptyset.
\eeq
\end{say}
Consequently, we have

\begin{prop} Fix any $m \ge d-1$. $N_d^l$ is the proper transform
of $N_d$ in $S^{m+1+l}(W_m,V_{d+m})$. In particular, $M_d$ is the
proper transform of $N_d$ in $S^{d+m}(W_m,V_{d+m})$
\end{prop}

Further, we have

\begin{lemm}
 Assume $1 \le l < k$. Then
$$I_{Z_{d,k}^{l-1}}\cdot
\sO_{N_d^l}=I_{Z_{d,k}^l}\cdot (I_{Z_{d,l}^l})^{k-l+1}.$$
\end{lemm}
\begin{proof} In the proof, we fix $m=d-1$ and use the embedding
$$N_d \lra  S(W_{d-1},V_{2d-1}).$$ We will use the notations
introduced in \ref{sEF} with $E=W_{d-1}$ and $F=V_{2d-1}$.

By \eqref{empty}, we have
 $$N_d\cap D_1=\cdots= N_d\cap D_{d-1}=\emptyset.$$ Thus, we have the induced embedding
\[
N_d \lra  S^{d-1},
\]
and moreover, for $0\leq k\leq d-1$, we have
\[
N_d\cap D_{d+k}^{d-1}=Z_{d,k}.
\]
In other words,
\[
I_{D_{d+k}^{d-1}}\cdot \sO_{N_d}=I_{Z_{d,k}}.
\]
Thus, we have the following blowing-up diagram
\[
\begin{CD}
N_d^{l} @>>>S^{d+l} \\
@VVV @VVV\\
N_d^{l-1} @>>> S^{d+l-1}
\end{CD}
\]
for all $0 \le l \le d-1$. Further, we have
$$N_d^l \cap D_{d+k}^{d+l}=Z_{d,k}^l, \quad \hbox{for $l \le k$}.$$
Thus for $l < k$,
\[
I_{Z_{d,k}^{l-1}}\cdot \sO_{N_d^l}=(I_{D_{d+k}^{d+l-1}}\cdot
\sO_{N_d^{l-1}})\cdot \sO_{N_d^l}=
(I_{D_{d+k}^{d+l-1}}\cdot\sO_{S^{d+l}})\cdot \sO_{N_d^l}.
\]
By Proposition \ref{vain}, we have $I_{D_{d+k}^{d+l-1}}
\cdot\sO_{S^{d+l}}= I_{D_{d+k}^{d+l}}\cdot
(I_{D_{d+l}^{d+l}})^{k-l+1}$. It then follows that
\[
I_{Z_{d,k}^{l-1}}\cdot \sO_{N_d^l}=I_{D_{d+k}^{d+l}}\cdot
(I_{D_{d+l}^{d+l}})^{k-l+1}\cdot \sO_{N_d^l}= I_{Z_{d,k}^l}\cdot
(I_{Z_{d,l}^l})^{k-l+1}.
\]
\end{proof}

Applying the above lemma repeatedly, we obtain

\begin{coro}\lab{keyCoro} For all $0 \le l \le d-1$,
$$I_{Z_{d,l+1}}\cdot
\sO_{N_d^l}=(I_{Z_{d,l+1}^l})\cdot(I_{Z_{d,l}^l})^2\cdots
(I_{Z_{d,1}^l})^{l+1}\cdot(I_{Z_{d,0}^l})^{l+2}
 $$
\end{coro}

We remark here that $Z_{d,t}^l$ are Cartier divisors when $t \le
l$. This implies that the blowup of $N_d^l$ along the proper
transform $Z_{d,l+1}^l$ of $Z_{d,l+1} \subset N_d$ is the same as
the blowup of $N_d^l$ along the total transform of $Z_{d,l+1}
\subset N_d$.

\subsection{Main theorems and proofs}

\begin{theo}\lab{main:thm1}
 Let $-1 \le k \le d-1$. Then
\begin{enumerate}
\item  $N_d^{k}$ is nonsingular;
\item $N_d^k$ is isomorphic to the closure of the graph of the
rational map
\[
\begin{split}
N_d &\dashrightarrow \prod_{l=0}^k S(\bigwedge^{m+2+l}
W_{m},\bigwedge^{m+2+l} V_{d+m}) \\
[f] &\mapsto
([\bigwedge^{m+2}\rho_{f,m}],[\bigwedge^{m+3}\rho_{f,m}],\dots,[\bigwedge^{m+2+k}\rho_{f,m}])
\end{split}
\]
for all $m \ge d-1$; \item $Z_{d,k+1}^k$ is isomorphic to
$\mathbb{P}(V_{d-k-1})\times N_{k+1}^k$.
\end{enumerate}
\end{theo}
\begin{proof}
We prove it by induction on $k$. When $k=-1$, the statements of the theorem
are trivial (for all $d>0$ and $m \ge d-1$).
Assume the statements are true for all $\leq k-1$ (for all $d>0$ and $m \ge d-1$).
We now prove the $k$-version of the theorem.

First, by the inductive assumption,  $N_{d}^{k-1}$ is nonsingular; also,
$Z_{d,k}^{k-1}$ is smooth because it is isomorphic to $ \mathbb{P}(V_{d-k-1})\times N_{k}^{k-1}$.
 Hence, $N_d^k$, as the blowup of  $N_d^{k-1}$ along  $Z_{d,k}^{k-1}$, is nonsingular.
 So (1) holds true for $k$.

Next, we let
 $$\pi_{[k-1]}: N_d^{k-1}\to N_d$$
be the iterated blowing-up morphism. By Corollary \ref{keyCoro},
the blowup of $N_d^{k-1}$ along the proper transform
$Z_{d,k}^{k-1}$ is isomorphic to
 the the blowup of $N_d^{k-1}$ along the total transform $\pi_{[k-1]}^{-1}(Z_{d,k})$.
Hence, by the definition of $Z_{d,k}$ (see Definition
\ref{defn:zdk} and Proposition \ref{weaker-sameIdeals}), $N_d^k$
is isomorphic to the closure of the graph of the rational map
\[
N_d^{k-1}\dashrightarrow S(\bigwedge^{k+2+m}
W_m,\bigwedge^{k+2+m}V_{d+m}).
\]
 By the induction hypothesis, $N_d^{k-1}$ is
isomorphic to the closure of the graph of the rational map
\[
N_d \dashrightarrow \prod_{l=0}^{k-1} S(\bigwedge^{m+2+l}
W_{m},\bigwedge^{m+2+l} V_{d+m})).
\]
It follows that $N_d^k$ is isomorphic to the closure of graph of
the rational map
\[
N_d \dashrightarrow \prod_{l=0}^{k} S(\bigwedge^{m+2+l}
W_{m},\bigwedge^{m+2+l} V_{d+m})).
\]
Thus (2) also holds true for $k$.

Finally, to prove the $k$-version  of (3), we  introduce and
establish the following commutative diagram
\beq\lab{commutativeDiagram}
\begin{CD}
Z_{d,k+1}\setminus Z_{d,k} @>{\alpha}>>\mathbb{P}(V_{d-k-1})\times
N_{k+1}^{k-1}\times
S(\ds\bigwedge^{k+m+2}W_m,\bigwedge^{k+m+2}V_{k+m+1}) \\
@| @VV{\beta}V\\
Z_{d,k+1}\setminus Z_{d,k} @>{\gamma}>> N_d^{k-1}\times
S(\ds\bigwedge^{k+m+2}W_m,\bigwedge^{k+m+2}V_{d+m}).
\end{CD}
\eeq Here $\gamma$ is the obvious embedding. The morphism $\alpha$
is the obvious embedding induced by the inverse of the morphism
$\varphi_{d,k+1}'$ of Proposition \ref{prop:inverse}. The morphism
$\beta$ is defined as follows.

For any integers $r,s\geq 1$, using the following identifications
\[\begin{split}& V_{r} = \Hom(\sO_{\Po}(-r),\sO_{\Po}),\\
& V_s= \Hom(\sO_{\Po},\sO_{\Po}(s)),\\
& V_{r+s}\Hom(\sO_{\Po}(-r),\sO_{\Po}(s)),
\end{split}\]
we see that each nonzero $h\in V_r$ gives rise to an injective
linear map
\[\begin{split}
& L_h: V_s \lra  V_{r+s} \\
& g \mapsto g\circ h.
\end{split}
\]
For any $1\leq l\leq s+1$, it induces an injective linear map
\[
\bigwedge^lL_h: \bigwedge^l V_s \to \bigwedge^l V_{r+s},
\]
 which in turn induces a morphism
\[
\mathbb{P}(V_r)\to S(\bigwedge^l V_s,\bigwedge^l V_{r+s}).
\]
Now fix any $t\geq 1$. Then, by the means of composing with
$\bigwedge^l L_h$, we obtain a morphism
\[
S(\bigwedge^l W_t,\bigwedge^l V_{s})\to S(\bigwedge^l
W_t,\bigwedge^l V_{r+s}).
\]
 Since $\bigwedge^l L_h$ is injective, the above morphism is an embedding (in fact, a linear embedding).

Observe now that when $l=s+1$, $\dim \bigwedge^l V_s=1$. Hence
$$S(\bigwedge^{s+1} V_s,\bigwedge^{s+1} V_{r+s})\cong
\mathbb{P}(\bigwedge^{s+1} V_{r+s}).$$  Then one checks directly
that the morphism
\[\mathbb{P}(V_r)\to S(\bigwedge^{s+1} V_s,\bigwedge^{s+1} V_{r+s})\] is the Veronese
embedding; further, the morphism
\[
S(\bigwedge^{s+1}W_t,\bigwedge^{s+1}V_s)\times
S(\bigwedge^{s+1}V_s,\bigwedge^{s+1}V_{r+s}) \lra
S(\bigwedge^{s+1}W_t,\bigwedge^{s+1}V_{r+s})
\]
is the  Segre embedding.

We are now ready to define $\beta$. It consists of two components
$(\beta_1, \beta_2)$. They are
\[\begin{split} & \beta_1: \mathbb{P}(V_{d-k-1})\times N_{k+1}^{k-1}\to
N_d^{k-1} \\
& \beta_1 ([h],[g],\prod_{l=0}^{k-1}[\bigwedge^{m+2+l}\rho_{g,m}])
= ([g \circ
h],\prod_{l=0}^{k-1}[\bigwedge^{m+2+l}L_h\circ\bigwedge^{m+2+l}\rho_{g,m}])
\end{split}
\]
and
\[\begin{split}
& \beta_2: \mathbb{P}(V_{d-k-1})\times S(\ds\bigwedge^{k+m+2}W_m,
\bigwedge^{k+m+2}V_{k+m+1}) \to
S(\ds\bigwedge^{k+m+2}W_m,\bigwedge^{k+m+2}V_{d+m})\\
& \beta_2([h],[\bigwedge^{m+2+k}\rho_{g,m}]) =
[\bigwedge^{m+2+l}L_h\circ\bigwedge^{m+2+k}\rho_{g,m}].
\end{split}
\]
Then one checks routinely that $\beta_2$ is the composition of the
Veronese embedding
\[\proj(V_{d-k+1}) \lra S(\bigwedge^{k+m+2}V_{k+m+1},
\bigwedge^{k+m+2}V_{d+m})\] followed by the Segre embedding
\[ \begin{split} & S(\ds\bigwedge^{k+m+2}W_m,
\bigwedge^{k+m+2}V_{k+m+1}) \times S(\bigwedge^{k+m+2}V_{k+m+1},
\bigwedge^{k+m+2}V_{d+m})  \\
& \lra S(\ds\bigwedge^{k+m+2}W_m,\bigwedge^{k+m+2}V_{d+m}),
\end{split}
\]
 hence
an embedding itself.  Since it is routine, we omit the details.
For $\beta_1$, observe that when $[h]$ is fixed,  it is injective
on other factors. Together, this implies that $\beta$ is an
embedding. Thus, we finally established the embedding diagram
\eqref{commutativeDiagram}.

From \eqref{commutativeDiagram}, we see that the closure of
$Z_{d,k+1}\setminus Z_{d,k}$ in
$$\mathbb{P}(V_{d-k-1})\times
N_{k+1}^{k-1}\times
S(\ds\bigwedge^{k+m+2}W_m,\bigwedge^{k+m+2}V_{k+m+1})$$  is
contained in $\mathbb{P}(V_{d-k-1})\times N_{k+1}^k$, and hence
equals to $\mathbb{P}(V_{d-k-1})\times N_{k+1}^k$ because
$Z_{d+1}\setminus Z_{d,k}$ is an open subvariety in it. Because
the embedding diagram \eqref{commutativeDiagram} commutes, the
above-mentioned closure is obviously isomorphic to the closure of
$Z_{d,k+1}\setminus Z_{d,k}$ in $$N_d^{k-1}\times
S(\ds\bigwedge^{k+m+2}W_m,\bigwedge^{k+m+2}V_{d+m}),$$ which is by
definition $Z_{d,k+1}^k$. This proves (3) for $k$.

By induction, the theorem is proved.
\end{proof}

\begin{theo}\lab{main:thm2}
Let $-1 \le k \le d-1$. Then

(1) $Z_{d,0}^k\cup\dots\cup Z_{d,k}^k$ is a divisor of $N_d^k$
with simple normal crossings;

(2) The scheme-theoretic intersection $Z_{d,k+1}^k\cap
\bigcap_{j=1}^r Z_{d,i_j}^k$ is isomorphic to
\[
\mathbb{P}(V_{d-k-1})\times\bigcap_{j=1}^r Z_{k+1,i_j}^k
\]
for any  distinct integers $i_1,\dots,i_r$ between $0$ and $k$.
\end{theo}
\proof Again, we prove it by induction on $k$. When $k=-1$, both
statements are trivial.

Assume that both statements are true for $k-1$. Then, it implies that $Z_{d,0}^{k-1}\cup \cdots
\cup Z_{d,k-1}^{k-1}$ is a divisor of $N_d^{k-1}$ with simple normal crossings. { In addition, by
statement (2) from the induction hypothesis, the scheme-theoretic intersection of $Z_{d,k}^{k-1}$
with an arbitrary intersection of $Z_{d,i_j}^{k-1}$'s is smooth. Next, it is routine to check that
\[\begin{split}
{\rm codim}\bigg(Z_{d,k}^{k-1}\cap\bigcap_{j=1}^r Z_{d,i_j}^{k-1}\bigg)&=\dim N_d^{k-1}-\dim\bigg(
\mathbb{P}(V_{d-k})\times \bigcap_{j=1}^r Z_{k,i_j}^{k-1}\bigg)\\ &={\rm
codim}Z_{d,k}^{k-1}+\sum_{j=1}^r{\rm codim} Z_{d,i_j}^{k-1}.
\end{split}
\]
Since it is routine, we omit further details.
Thus,  $Z_{d,k}^{k-1}$ meets the divisors $Z_{d,0}^{k-1}, \ldots, Z_{d,k-1}^{k-1}$
transversally.} Since transversality is preserved under blowup along nonsingular center, we obtain
the $k$-version of (1).

To prove (2), we consider again the commutative diagram in
\eqref{commutativeDiagram}, which by now can also be expressed as
$$
\begin{CD}
Z_{d,k+1}^k @>{\alpha}>>\mathbb{P}(V_{d-k-1})\times
N_{k+1}^{k-1}\times
S(\ds\bigwedge^{k+m+2}W_m,\bigwedge^{k+m+2}V_{k+m+1}) \\
@| @VV{\beta}V\\
Z_{d,k+1}^k  @>{\gamma}>> N_d^{k-1}\times
S(\ds\bigwedge^{k+m+2}W_m,\bigwedge^{k+m+2}V_{d+m}).
\end{CD}
$$
Then  from the diagram we have that for any $0 \le i \le k$,
$$Z_{d,k+1}^k \cap Z_{d,i}^k
= \proj(V_{d-k-1}) \times Z_{k+1,i}^{k}$$ scheme-theoretically and
further the scheme theoretic intersection
$$Z_{d,k+1}^k \cap \bigcap_{j=1}^r Z_{d,i_j}^k$$ is
isomorphic to
$$\proj(V_{d-k-1}) \times \bigcap_{j=1}^r Z_{k+1,i_j}^k.$$
\endproof

Theorems \ref{main:thm1} and \ref{main:thm2} imply

\begin{theo}\lab{main:grand} Fix any $m \ge d-1$. Then
$M_d = N_d^{d-1}$ is a  compactification of $\mathring{N}_d$ such
that the following hold.
\begin{enumerate}
\item $M_d$ is isomorphic to the closure of the graph of the
rational map
\[
N_d \dashrightarrow \prod_{l=0}^{d-1} S(\bigwedge^{m+2+l}
W_{m},\bigwedge^{m+2+l} V_{d+m});
\]
\item $M_d$ is a nonsingular projective variety; \item The
complement $M_d \setminus \mathring{N}_d  = \bigcup_{k=0}^{d-1}
Z_{d,k}^{d-1}$ is a divisor with simple normal crossings.
\end{enumerate}
\end{theo}

\section{The Topology of the Compactification}

\begin{say} In this section, we will work over the field of complex numbers  and draw a
few consequences on the topology of $M_d$. Recall that for any
complex quasi-projective variety $V$ there is a (virtual Hodge)
polynomial $\ee (V)$ in two variables $u$ and $v$ which is
uniquely determined by the following properties
\begin{enumerate}
\item If $V$ is smooth and projective, then $\ee (V) = \sum h^{p,
q} (-u)^p (-v)^q$; \item If $U$ is a closed subvariety of $V$,
then $\ee (V) = \ee (V     \setminus  U) + \ee (U)$; \item If $V
\rightarrow B$ is a Zariski locally trivial bundle with fiber $F$,
then $\ee (V) = \ee (B) \ee (F)$.
\end{enumerate}
\end{say}

\begin{say}
Fix $n>0$.  For any $i >0$, set
$$R_i (\lambda) = \frac{\lambda^{i+1} -1}{\lambda-1} \cdot
\frac{\lambda^{ni} -\lambda}{\lambda-1}.$$ First, we have the
following recursive formula for $\ee_{M_d}$.
\end{say}

\begin{prop}
\label{recursiveformula} Set $\lambda = uv$. Then
$$\ee_{M_d} = \ee_{N_d} + \sum_{k=0}^{d-1}
\ee_{M_k} R_{d-k}$$ where $\ee_{N_d} = \frac{\lambda^{(d+1)(n+1)}
-1}{\lambda-1}$.
\end{prop}
\proof Since $M_d =N_d^{d-1}$ is the blowup of $N_d^{d-2}$ along
$Z_{d,d-1}^{d-2}$, we have
$$\ee_{M_d} = \ee_{N_d^{d-2}} +
\ee_{Z_{d,d-1}^{d-2}} (\ee_{\proj^{(\codim Z_{d,d-1}-1)}} -1),$$
that is
$$\ee_{M_d} = \ee_{N_d^{d-2}} +
\ee_{Z_{d,d-1}^{d-2}} \frac{\lambda^{\codim Z_{d,d-1}} -
\lambda}{\lambda -1}.$$ Repeat the same arguments for $\ee_{
N_d^{d-2}}$ and so on, we will eventually obtain
$$\ee_{M_d} =  \ee_{N_d} + \sum_{k=0}^{d-1}
\ee_{Z_{d,k}^{k-1}} \frac{\lambda^{\codim Z_{d,k}} -
\lambda}{\lambda -1}.$$ Because $Z_{d,k}^{k-1} = \proj(V_{d-k})
\times M_k$, $\proj(V_{d-k}) \cong \proj^{d-k}$ and $\codim
Z_{d,k} = n(d-k)$, from here it is routine to  obtain the formula
as stated in the theorem.
\endproof

\begin{say} We can also derive a closed formula for $\ee_{M_d}$.
For this
 consider $\alpha= (\alpha_1, \ldots,
\alpha_r) \in \ZZ_{>0}^r$. Let $|\alpha|$ denote the sum
$\sum_{j=1}^r \alpha_j$.  Set
$$R_\alpha = \prod_{j=1}^r R_{\alpha_j} \and R_0 = R_{0,n} =1.$$
\end{say}

\begin{prop}
\label{closeformula} The Hodge  polynomial $\ee_{M_d}(u,v)$ is
given by
$$\ee_{M_d} = \sum_{ 0 \le |\alpha| \le d}
 R_\alpha \ee_{N_{d-|\alpha|}} $$ where $\ee_{N_k} =
\frac{\lambda^{(k+1)(n+1)} -1}{\lambda-1}$ for all $k \ge 0$.
\end{prop}
\proof When $d=0$, it is trivial. Assume that the formula holds
for all $k < d$. By Theorem \ref{recursiveformula},
$$\ee_{M_d} = \ee_{N_d} + \sum_{k=0}^{d-1}
\ee_{M_k} R_{d-k}.$$ By inductive assumption, for all $k < d$
$$\ee_{M_k} = \sum_{ 0 \le |\alpha| \le k}
 R_\alpha \ee_{N_{k-|\alpha|}} 
 = \sum_{\scriptsize {\begin{array}{cccl}  0 \le |\alpha| \le k \\
  |\alpha| +j =k \end{array}}}
 R_\alpha \ee_{N_j} .$$ Substitute $\ee_{M_k}$ into the first
formula, we have
$$
\ee_{M_d}  = \ee_{N_d} + \sum_{k=0}^{d-1}
\sum_{\scriptsize{ \begin{array}{cccl}  0 \le |\alpha| \le k \\
  |\alpha| +j =k \end{array}}}
 R_\alpha  R_{d-k} \ee_{N_j}.$$
Let $\beta= (\alpha, d-k)$. Then $R_\beta = R_\alpha R_{d-k}$ and
$\beta =|\alpha| + d-k = d-j$ if $|\alpha| = k-j$. Then a simple
counting argument shows
$$\ee_{M_d} = \ee_{N_d} +
\sum_{\scriptsize{ \begin{array}{cccl} |\beta| + j= d \\
0 \le j \le d-1 \end{array}}}  R_\beta \ee_{N_j} 
= \sum_{\scriptsize{ \begin{array}{cccl}  0 \le |\beta| \le d \\
  |\beta| +j =d \end{array}}}
 R_\beta \ee_{N_j} = \sum_{0 \le |\beta| \le d} R_\beta \ee_{N_{d-|\beta|}} ,$$ as desired.
\endproof

In addition, we have

\begin{prop}
Let $H$ be the pullback (to $M_d$) of the hyperplane class of
$N_d$ and $T_k$ corresponds to the exceptional divisor
$Z_{d,k}^{d-1}$ for all $0 \le k \le d-1$. Then $$A^1 (M_d) = \ZZ
\cdot H \oplus \bigoplus_{k=0}^{d-1} \ZZ \cdot T_k.$$
\end{prop}
\proof This follows from that fact that $M_d$ is a successive
blowup of $N_d$ along nonsingular centers.
\endproof

\end{document}